\documentclass[12pt]{article}

\usepackage{amsmath, epsfig, cite}
\usepackage{amssymb}
\usepackage{amsfonts}
\usepackage{latexsym}
\usepackage{float}

\newtheorem{thm}{Theorem}[section]

\newtheorem{conj}[thm]{Conjecture}
\newtheorem{lem}[thm]{Lemma}

\numberwithin{equation}{section}

\newcommand{\qed}{{\hfill$\square$}\medskip}

\begin{document}


\begin{center}
{\Large\bf  A supercongruence involving Delannoy numbers\\[10pt]
and Schr\"oder numbers}
\end{center}

\vskip 2mm \centerline{Ji-Cai Liu}
\begin{center}
{\footnotesize Department of Mathematics, Shanghai Key Laboratory of
PMMP, East China Normal University,\\ 500 Dongchuan Road, Shanghai
200241,
 People's Republic of China\\
{\tt jc2051@163.com} }
\end{center}


\vskip 0.7cm \noindent{\bf Abstract.}
The Delannoy numbers and Schr\"oder numbers are given by
\begin{align*}
D_n=\sum_{k=0}^n{n\choose k}{n+k\choose k}\quad \text{and}\quad
S_n=\sum_{k=0}^n{n\choose k}{n+k\choose k}\frac{1}{k+1},
\end{align*}
respectively. Let $p>3$ be a prime. We mainly prove that
\begin{align*}
\sum_{k=1}^{p-1}D_k S_k\equiv 2p^3B_{p-3}-2pH^{*}_{p-1} \pmod{p^4},
\end{align*}
where $B_n$ is the $n$-th Bernoulli number and those $H^{*}_n$ are the alternating harmonic numbers given by $H^{*}_n=\sum_{k=1}^{n}\frac{(-1)^k}{k}$.
This supercongruence was originally conjectured by Z.-W. Sun in 2011.
\vskip 2mm
\noindent{\it MR Subject Classifications}: 11A07, 11B65, 11B83, 05A10

\section{Introduction}
In combinatorics, the $n$-th Delannoy number describes the number of paths from $(0,0)$ to $(n,n)$, using only steps $(1,0), (0,1)$ and $(1,1)$, while $n$-th Schr\"oder number represents the number of such paths that do not rise above the line $y=x$. It is known that
\begin{align*}
D_n=\sum_{k=0}^{n}{n\choose k}{n+k\choose k}\quad \text{and} \quad
S_n=\sum_{k=0}^{n}{n\choose k}{n+k\choose k}\frac{1}{k+1}.
\end{align*}

Z.-W. Sun\cite{Sun,Sun2} have proved some amazing arithmetic properties of Delannoy numbers and Schr\"oder numbers. For example, he showed that for any odd prime $p$,
\begin{align*}
&\sum_{k=1}^{p-1}\frac{D_k}{k^2}\equiv 2\left(\frac{-1}{p}\right)E_{p-3} \pmod{p},\\
&\sum_{k=0}^{p-1}D^2_k\equiv \left(\frac{2}{p}\right)\pmod{p},
\end{align*}
where $\left(\frac{\cdot}{p}\right)$ denotes the Legendre symbol and $E_n$ is the $n$-th Euler number. In 2011, Z.-W. Sun\cite{Sun} also raised some interesting conjectures involving those numbers, one of which was
\begin{conj}
Let $p>3$ be a prime. Then
\begin{align}
\sum_{k=1}^{p-1}D_k S_k\equiv -2p\sum_{k=1}^{p-1}\frac{(-1)^k+3}{k} \pmod{p^4}. \label{a3}
\end{align}
\end{conj}

Define the following two polynomials:
\begin{align*}
D_n(x)=\sum_{k=0}^{n}{n\choose k}{n+k\choose k}x^k\quad \text{and} \quad S_n(x)=\sum_{k=0}^{n}{n\choose k}{n+k\choose k}\frac{x^k}{k+1}.
\end{align*}
Note that $D_n(\frac{x-1}{2})$ coincides with the Legendre polynomial $P_n(x)$.
Recently, Guo\cite{Guo} proved that for any prime $p$ and integer $x$ satisfying $p\nmid x(x+1)$,
\begin{align*}
&\sum_{k=0}^{p-1}(2k+1)D_k(x)^3\equiv p\left(\frac{-4x-3}{p}\right) \pmod{p^2},\\
&\sum_{k=0}^{p-1}(2k+1)D_k(x)^4\equiv p \pmod{p^2},
\end{align*}
and thus confirmed some conjectures of Z.-W. Sun \cite{Sun2}.

Congruences for similar numbers and polynomials have been widely studied by several
authors, see, for example, Ahlgren and Ono \cite{AK}, Gessel\cite{Ge} and Pan\cite{Pan}.

Our first result is the following congruence:
\begin{thm}\label{t1}
Let $p$ be an odd prime and $x$ be an integer not divisible by $p$. Then
\begin{align}
\sum_{k=1}^{p-1}D_k(x)S_k(x)\equiv 0 \pmod{p}.\label{a1}
\end{align}
\end{thm}

Note that, by the equation below (16) in Lehmer's paper\cite{Leh},
\begin{align}
\sum_{k=1}^{p-1}\frac{1}{k}\equiv -\frac{p^2}{3}B_{p-3} \pmod{p^3}, \label{d0}
\end{align}
where $B_n$ is the $n$-th Bernoulli number given by
\begin{align*}
\frac{x}{e^x-1}=\sum_{n=0}^{\infty}B_n\frac{x^n}{n!},\quad \text{for $0< |x|< 2\pi $}.
\end{align*}
From \eqref{d0}, it is clear that \eqref{a3} is equivalent to the following supercongruence:
\begin{thm}\label{t2}
Let $p>3$ be a prime. Then
\begin{align}
\sum_{k=1}^{p-1}D_k S_k\equiv 2p^3B_{p-3}-2pH^{*}_{p-1} \pmod{p^4}, \label{a2}
\end{align}
where $H^{*}_n=\sum_{k=1}^{n}\frac{(-1)^k}{k}$.
\end{thm}
We will prove this supercongruence in Section 3.

\section{Proof of Theorem \ref{t1}}

We begin with the following combinatorial identity (see \cite[(2.5)]{Guo}):
\begin{align}
{k\choose i}{k+i\choose i}{k\choose j}{k+j\choose j}
&=\sum_{r=0}^i {i+j\choose i}{j\choose i-r}{j+r\choose r}{k\choose j+r}{k+j+r\choose j+r} \notag\\
&=\sum_{s=\text{max}\{i,j\}}^{i+j} {s\choose i}{s\choose j}{i+j\choose s}{k\choose s}{k+s\choose s}.\label{b1}
\end{align}
Using \eqref{b1} and exchanging summation order, we have
\begin{align}
&\hskip -2mm
\sum_{k=0}^{p-1}D_k(x)S_k(x)\notag\\
&=\sum_{k=0}^{p-1}\sum_{i=0}^k\sum_{j=0}^k{k\choose i}{k+i\choose i}{k\choose j}{k+j\choose j}\frac{x^{i+j}}{j+1}\notag\\
&=\sum_{k=0}^{p-1}\sum_{i=0}^k\sum_{j=0}^k\sum_{s} {s\choose i}{s\choose j}{i+j\choose s}{k\choose s}{k+s\choose s}\frac{x^{i+j}}{j+1}\notag\\
&=\sum_{i=0}^{p-1}\sum_{j=0}^{p-1}\sum_{s} {s\choose i}{s\choose j}{i+j\choose s}\frac{x^{i+j}}{j+1}\sum_{k=\text{max}\{i,j\}}^{p-1}{k+s\choose s}{k\choose s}.\label{b2}
\end{align}
Noticing that $s\ge \text{max}\{i,j\}$ in \eqref{b2} and applying the following identity:
\begin{align*}
\sum_{k=s}^{n-1}{k\choose s}{k+s\choose s}=\frac{n}{2s+1}{n+s\choose s}{n-1\choose s},
\end{align*}
which can be easily proved by induction on $n$, we get
\begin{align}
&\hskip -2mm
\sum_{k=0}^{p-1}D_k(x)S_k(x)\notag\\
&=p\sum_{i=0}^{p-1}\sum_{j=0}^{p-1}\sum_{s=0}^{p-1}  {s\choose i}{s\choose j}{i+j\choose s}{p+s\choose s}{p-1\choose s}\frac{x^{i+j}}{(j+1)(2s+1)}\notag\\
&=p\sum_{s=0}^{p-1}\sum_{m=0}^{2p-2}\sum_{j=0}^{p-1} {s+1\choose j+1}{s\choose m-j}{m\choose s}{p+s\choose s}{p-1\choose s}\frac{x^{m}}{(s+1)(2s+1)}, \label{b3}
\end{align}
where $m=i+j$.
By Chu-Vandermonde formula, we have
\begin{align}
\sum_{j=0}^{s} {s+1\choose j+1}{s\choose m-j}={2s+1\choose m+1},\quad \text{for $s\le m$}.\label{f1}
\end{align}
It follows from \eqref{b3} and \eqref{f1} that
\begin{align}
&\hskip -2mm
\sum_{k=0}^{p-1}D_k(x)S_k(x)\notag\\
&=p\sum_{s=0}^{p-1}\sum_{m=0}^{2p-2}{2s+1\choose m+1}{m\choose s}{p+s\choose s}{p-1\choose s}\frac{x^m}{(s+1)(2s+1)}\notag\\
&=p\sum_{m=0}^{2p-2}\sum_{s=0}^{p-1}{2s\choose m}{m+1\choose s+1}{p+s\choose s}{p-1\choose s}\frac{x^m}{(m+1)^2} \label{b4}.
\end{align}
Since the summands in right-hand side of \eqref{b4} are congruent to $0$ mod $p$ for $m\neq p-1$, we immediately get
\begin{align*}
\sum_{k=0}^{p-1}D_k(x)S_k(x)
\equiv \frac{x^{p-1}}{p}\sum_{s=0}^{p-1}{2s\choose p-1}{p\choose s+1}{p+s\choose s}{p-1\choose s} \pmod{p}.
\end{align*}
Note that $\frac{1}{p}{2s\choose p-1}{p\choose s+1}$ is always $p$-adically integral for $0\le s\le p-1$ and ${p+s\choose s}{p-1\choose s}\equiv (-1)^s \pmod{p}$. Thus,
\begin{align*}
\sum_{k=0}^{p-1}D_k(x)S_k(x)
&\equiv \frac{x^{p-1}}{p}\sum_{s=0}^{p-1}(-1)^s{2s\choose p-1}{p\choose s+1} \pmod{p}\\
&=x^{p-1}\sum_{s=0}^{p-1}{2s\choose s}{s\choose p-1-s}\frac{(-1)^s}{s+1}.
\end{align*}
Furthermore, applying the following identity
\begin{align}
\sum_{s=0}^{m}{2s\choose s}{s\choose m-s}\frac{(-1)^s}{s+1}=(-1)^m, \label{b0}
\end{align}
which can be proved by Zeilberger's algorithm (see \cite{PWZ}),
and then noting that $D_0(x)S_0(x)=1$ and $x^{p-1} \equiv 1 \pmod{p}$ for $p\nmid x$,
we complete the proof.

\section{Proof of Theorem \ref{t2}}

Let
$$H^{(2)}_n=\sum_{k=1}^n\frac{1}{k^2}$$
be the $n$-th generalized harmonic number of order $2$.
The following lemmas will be used in our proof of Theorem \ref{t2}.
\begin{lem}
Let $m$ be a nonnegative integer. Then
\begin{align}
\sum_{s=0}^{m}{2s\choose s}{s\choose m-s}\frac{(-1)^s }{s+1}H^{(2)}_s=\frac{2(-1)^m}{m+1}\sum_{s=0}^m H^{(2)}_s. \label{c1}
\end{align}
\end{lem}
\noindent{\it Proof.}
Denote both sides of \eqref{c1} by $S_1(m)$ and $S_2(m)$, respectively. Applying the multi-variable Zeilberger algorithm (see \cite{AZ}), we find that
$S_1(m)$ and $S_2(m)$ satisfy the same recurrence:
\begin{align*}
&(m+4)(m+3)^2S_{i}(m+3)+(m+3)(3m^2+16m+22)S_{i}(m+2)\\
&+(m+2)(3m^2+14m+17)S_{i}(m+1)+(m+1)(m+2)^2S_{i}(m)=0, \quad \text{for $i=1,2$.}
\end{align*}
Also, it is easily verified that $S_1(m)=S_2(m)$ for $m=1,2,3$. This proves the lemma.
\qed

\begin{lem}
Let $p>3$ be a prime and $m$ be a positive integer such that $p\le m \le 2p-2$. Then
\begin{align}
&\sum_{s=0}^{p-1}{2s\choose s}{s\choose m-s}\frac{(-1)^s}{s+1}\notag\\
&\equiv (-1)^m\left(-1+\frac{2p}{m+1}+\frac{2p^2}{m+1}\sum_{s=0}^{m-p}H^{(2)}_s \right) \pmod{p^3}.\label{c2}
\end{align}
\end{lem}
\noindent{\it Proof.}
By \eqref{b0}, we have
\begin{align}
\sum_{s=0}^{p-1}{2s\choose s}{s\choose m-s}\frac{(-1)^s}{s+1}
=(-1)^{m}-\sum_{s=p}^{m}{2s\choose s}{s\choose m-s}\frac{(-1)^s}{s+1}.\label{e1}
\end{align}
Letting $m=p+r$, $0\le r\le p-2$ and $s\to s+p$ in the right-hand side of \eqref{e1}, we find that \eqref{c2} is equivalent to
\begin{align}
&\sum_{s=0}^r{2p+2s\choose p+s}{p+s\choose r-s}\frac{(-1)^s}{p+s+1}\notag\\
&\equiv 2(-1)^r\left(1-\frac{p}{p+r+1}-\frac{p^2}{p+r+1}\sum_{s=0}^r H^{(2)}_s\right) \pmod{p^3}.\label{e0}
\end{align}
Note that
\begin{align}
{2p+2s\choose p+s}{p+s\choose r-s}&=\frac{2^s}{(r-s)!}{2p\choose p}\prod_{i=1}^s(2p+2i-1)
\left. \middle/\right.\prod_{i=1}^{2s-r}(p+i)\notag\\
&\equiv \frac{2^{s+1}}{(r-s)!}\prod_{i=1}^s(2p+2i-1)\left. /\right.\prod_{i=1}^{2s-r}(p+i) \pmod{p^3},\label{e2}
\end{align}
where we have used Wolstenholme's theorem ${2p\choose p}\equiv 2\pmod{p^3}$ for $p\ge 5$.
Now consider the following rational function:
\begin{align*}
f_{r,s}(x)=\frac{(-2)^s}{(r-s)!(x+s+1)}\prod_{i=1}^s(2x+2i-1)\left./\right.\prod_{i=1}^{2s-r}(x+i).
\end{align*}
Taking the logarithmic derivative on both sides, we obtain
\begin{align}
f'_{r,s}(x)
=f_{r,s}(x)\left(\sum_{i=1}^s\frac{2}{2x+2i-1}-\sum_{i=1}^{2s-r}\frac{1}{x+i}-\frac{1}{x+s+1}\right).
\label{e3}
\end{align}
Differentiating both sides of \eqref{e3}, we get the second order derivative:
\begin{align*}
f''_{r,s}(x)&=f_{r,s}(x)\left(\sum_{i=1}^s\frac{2}{2x+2i-1}-\sum_{i=1}^{2s-r}\frac{1}{x+i}-\frac{1}{x+s+1}\right)^2\\
&+f_{r,s}(x)\left(\frac{1}{(x+s+1)^2}+\sum_{i=1}^{2s-r}\frac{1}{(x+i)^2}-\sum_{i=1}^{s}\frac{4}{(2x+2i-1)^2}\right).
\end{align*}
It is not hard to see that
\begin{align}
&f_{r,s}(0)={2s\choose s}{s\choose r-s}\frac{(-1)^s}{s+1},\notag\\
&\sum_{s=0}^r f_{r,s}(0)=(-1)^r. \quad \text{(by \eqref{b0})}\label{e4}
\end{align}
Similarly to the proof of \eqref{c1}, we can obtain the following two identities:
\begin{align}
&\sum_{s=0}^r f'_{r,s}(0)=\frac{(-1)^{r+1}}{r+1},\label{e5}\\
&\sum_{s=0}^r f''_{r,s}(0)=\frac{2(-1)^{r}}{(r+1)^2}+\frac{2(-1)^{r+1}}{r+1}\sum_{s=0}^r H^{(2)}_s.\label{e6}
\end{align}
Now combining \eqref{e4}-\eqref{e6}, we get the first three terms of Taylor expansion:
\begin{align}
\sum_{s=0}^rf_{r,s}(x)&=(-1)^r\left(1-\frac{x}{r+1}
+\frac{x^2}{(r+1)^2}-\frac{x^2}{r+1}\sum_{s=0}^r H^{(2)}_s\right)+\mathcal O\left(x^{3}\right)\notag\\
&=(-1)^r\left(1-\frac{x}{x+r+1}
-\frac{x^2}{x+r+1}\sum_{s=0}^r H^{(2)}_s\right)+\mathcal O\left(x^{3}\right).\label{e7}
\end{align}
The proof of \eqref{e0} directly follows from \eqref{e2} and \eqref{e7}.
\qed

\begin{lem}
Let $p>3$ be a prime and $m$ be a nonnegative integer. Then
\begin{align}
\sum_{s=0}^m H^{(2)}_s\equiv \sum_{s=0}^{p-m-2} H^{(2)}_s \pmod{p},\quad \text{for $0\le m \le p-2$}.\label{c3}
\end{align}
\end{lem}
\noindent{\it Proof.}
Note that
\begin{align}
H^{(2)}_s+H^{(2)}_{p-s-1}
&=\sum_{k=1}^{s}\frac{1}{k^2}+\sum_{k=s+1}^{p-1}\frac{1}{(p-k)^2}\notag\\
&\equiv \sum_{k=1}^{p-1}\frac{1}{k^2}\notag\\
&\equiv 0 \pmod{p}.\label{c3-1}
\end{align}
In particular, letting $s=\frac{p-1}{2}$ in \eqref{c3-1} yields $H^{(2)}_{\frac{p-1}{2}}\equiv 0\pmod{p}$.
Assuming $m\le \frac{p-3}{2}$ and using \eqref{c3-1}, we obtain
\begin{align*}
\sum_{s=0}^{p-m-2}H^{(2)}_s-\sum_{s=0}^{m}H^{(2)}_s
&=\sum_{s=m+1}^{\frac{p-1}{2}}\left(H^{(2)}_s+H^{(2)}_{p-s-1}\right)-H^{(2)}_{\frac{p-1}{2}}\\
&\equiv 0 \pmod{p}.
\end{align*}
This completes the proof.
\qed

\noindent{\it Proof of \eqref{a2}.}
Letting $x=1$ in \eqref{b4} reduces to
\begin{align}
\sum_{k=0}^{p-1}D_kS_k
&=p\sum_{m=0}^{2p-2}\sum_{s=0}^{p-1}{2s\choose m}{m+1\choose s+1}{p+s\choose s}{p-1\choose s}\frac{1}{(m+1)^2}\notag\\
&=p\sum_{m=0}^{2p-2}\sum_{s=0}^{p-1}{2s\choose s}{s\choose m-s}{p+s\choose s}{p-1\choose s}\frac{1}{(m+1)(s+1)}.\label{c4}
\end{align}
Note that for $0\le s\le p-1$,
\begin{align}
{p+s\choose s}{p-1\choose s}&=\frac{(p^2-s^2)\cdots (p^2-1^2)}{(s!)^2}\notag\\
&\equiv (-1)^s\left(1-H^{(2)}_s p^2\right) \pmod{p^4}. \label{c5}
\end{align}
Since for $0\le s\le p-1$ and $0\le m \le 2p-2$, $\frac{p}{(m+1)(s+1)}{2s\choose s}$ is always a $p$-adic integer, it follows from \eqref{c4} and \eqref{c5} that
\begin{align}
&\sum_{k=0}^{p-1}D_kS_k\notag\\
&\equiv p\sum_{m=0}^{2p-2}\sum_{s=0}^{p-1}\frac{(-1)^s{2s\choose s}{s\choose m-s}}{(m+1)(s+1)}
-p^3\sum_{m=0}^{2p-2}\sum_{s=0}^{p-1}\frac{(-1)^s{2s\choose s}{s\choose m-s}}{(m+1)(s+1)}H^{(2)}_s \pmod{p^4}. \label{c6}
\end{align}
We shall prove \eqref{a2} by establishing the following three supercongruences:
\begin{align}
&p\sum_{m=0}^{p-1}\sum_{s=0}^{p-1}\frac{(-1)^s{2s\choose s}{s\choose m-s}}{(m+1)(s+1)}
\equiv 1-pH^{*}_{p-1}\pmod{p^4},\label{k1}\\
&p\sum_{m=p}^{2p-2}\sum_{s=0}^{p-1}\frac{(-1)^s{2s\choose s}{s\choose m-s}}{(m+1)(s+1)}
\equiv 4p^3B_{p-3}-pH^{*}_{p-1} \pmod{p^4},\label{k2}\\
&p^3\sum_{m=0}^{2p-2}\sum_{s=0}^{p-1}\frac{(-1)^s{2s\choose s}{s\choose m-s}}{(m+1)(s+1)}H^{(2)}_s
\equiv 2p^3B_{p-3} \pmod{p^4}.\label{k3}
\end{align}
Substituting \eqref{k1}-\eqref{k3} into \eqref{c6} and noting that $D_0S_0=1$, we are led to Theorem \ref{t2}.
So it suffices to prove \eqref{k1}-\eqref{k3}.

By \eqref{b0}, we have
\begin{align*}
\text{LHS \eqref{k1}}
&=p\sum_{m=0}^{p-1}\frac{1}{m+1}\sum_{s=0}^{m}{2s\choose s}{s\choose m-s}\frac{(-1)^s}{s+1}\\
&=1-pH^{*}_{p-1},
\end{align*}
which proves \eqref{k1}.

Applying \eqref{c2} to the left-hand side of \eqref{k2} and then letting $m\to k+p-1$ yields
\begin{align*}
\text{LHS \eqref{k2}}\equiv p\sum_{k=1}^{p-1}\frac{(-1)^k}{k+p}\left(-1+\frac{2p}{k+p}+\frac{2p^2}{k+p}\sum_{s=0}^{k-1}H^{(2)}_s\right)
\pmod{p^4}.
\end{align*}
Note that
\begin{align*}
\frac{1}{k+p}\equiv \frac{1}{k}-\frac{p}{k^2}+\frac{p^2}{k^3} \pmod{p^3},\\
\frac{1}{(k+p)^2}=\frac{1}{k^2}-\frac{2p}{k^3} \pmod{p^2},
\end{align*}
and
\begin{align*}
\frac{(-1)^k}{(k+p)^2}\sum_{s=0}^{k-1}H^{(2)}_s+\frac{(-1)^{p-k}}{(2p-k)^2}\sum_{s=0}^{p-k-1}H^{(2)}_s
\equiv 0\pmod{p}. \quad \text{ (by \eqref{c3})}
\end{align*}
Hence
\begin{align}
\text{LHS \eqref{k2}}
\equiv -p\sum_{k=1}^{p-1}\frac{(-1)^k}{k}+3p^2\sum_{k=1}^{p-1}\frac{(-1)^k}{k^2}
-5p^3\sum_{k=1}^{p-1}\frac{(-1)^k}{k^3} \pmod{p^4}.\label{c7}
\end{align}
Furthermore, note that
\begin{align}
p^2\sum_{k=1}^{p-1}\frac{(-1)^k}{k^2}&=\frac{p^2}{2}\sum_{k=1}^{p-1}\left( \frac{(-1)^k}{k^2}+\frac{(-1)^{p-k}}{(p-k)^2}\right)\notag\\
&=\frac{p^3}{2}\sum_{k=1}^{p-1}\frac{(-1)^k(p-2k)}{k^2(p-k)^2}\notag\\
&\equiv -p^3\sum_{k=1}^{p-1}\frac{(-1)^k}{k^3} \pmod{p^4},\label{c8}
\end{align}
and
\begin{align}
\sum_{k=1}^{p-1}\frac{(-1)^k}{k^2}=\frac{1}{2}\sum_{k=1}^{\frac{p-1}{2}}\frac{1}{k^2}
-\sum_{k=1}^{p-1}\frac{1}{k^2}.\label{c9}
\end{align}
Combining \eqref{c7}-\eqref{c9}, we get
\begin{align}
\text{LHS \eqref{k2}}
\equiv -pH^{*}_{p-1}+4p^2\sum_{k=1}^{\frac{p-1}{2}}\frac{1}{k^2}-8p^2\sum_{k=1}^{p-1}\frac{1}{k^2}
\pmod{p^4}.\label{c10}
\end{align}
Now we go way back to Lehmer's paper \cite{Leh}. From the equations (15) and (18) in her paper, we can deduce the following congruence:
\begin{align*}
\sum_{k=1}^{\frac{p-1}{2}}k^{2r}\equiv \frac{1-2^{2r-1}}{2^{2r}}\sum_{k=1}^{p-1}k^{2r} \pmod{p^3}, \quad \text{for $2r\not\equiv 2 \pmod{p-1}$}.
\end{align*}
Letting $2r=p^3-p^2-2$ and using Euler's theorem to replace $x^{p^3-p^2-2}$ with $x^{-2}$, we get
\begin{align}
\sum_{k=1}^{\frac{p-1}{2}}\frac{1}{k^2}\equiv \frac{7}{2}\sum_{k=1}^{p-1}\frac{1}{k^2} \pmod{p^3},\quad \text{for $p\ge 7$}.\label{c11}
\end{align}
It is easily verified that the above congruence modulo $p^2$ is true for $p=5$.
Moreover, the following congruence can be deduced from the equation blow (16) in her paper:
\begin{align}
\sum_{k=1}^{p-1}\frac{1}{k^2}\equiv \frac{2p}{3}B_{p-3} \pmod{p^2}.\label{c12}
\end{align}
Substituting
\eqref{c11} into \eqref{c10} and then using \eqref{c12}, we complete the proof of \eqref{k2}.

Note that $H^{(2)}_{p-1}\equiv 0 \pmod{p}$ and ${2s\choose s}\equiv 0\pmod{p}$ for $\frac{p+1}{2}\le s\le p-1$. Thus, for $0\le s\le p-1$ and $p\le m\le 2p-2$, we have
\begin{align}
\frac{1}{s+1}{2s\choose s}{s\choose m-s}H^{(2)}_s\equiv 0 \pmod{p}.\label{f2}
\end{align}
Applying \eqref{f2} to the left-hand side of \eqref{k3} yields
\begin{align}
\text{LHS \eqref{k3}}
\equiv p^3\sum_{m=0}^{p-1}\frac{1}{m+1}\sum_{s=0}^{m}\frac{(-1)^s{2s\choose s}{s\choose m-s}}{s+1}H^{(2)}_s\pmod{p^4}.\label{c13}
\end{align}
Furthermore, substituting \eqref{c1} into the right-hand side of \eqref{c13} and noting that
for $0\le m\le p-2$,
\begin{align*}
\frac{(-1)^m}{(m+1)^2}\sum_{s=0}^m H^{(2)}_s +\frac{(-1)^{p-2-m}}{(p-1-m)^2}\sum_{s=0}^{p-2-m} H^{(2)}_s
\equiv 0\pmod{p}, \quad \text{(by \eqref{c3})}
\end{align*}
we immediately get
\begin{align*}
\text{LHS \eqref{k3}}
\equiv 2p\sum_{s=0}^{p-1} H^{(2)}_s \pmod{p^4}.
\end{align*}
Finally, noting that
\begin{align*}
\sum_{s=0}^{p-1} H^{(2)}_s=\sum_{k=1}^{p-1}\frac{p-k}{k^2}
\end{align*}
and then using \eqref{d0} and \eqref{c12},
we complete the proof of \eqref{k3}.
\qed

\vskip 5mm \noindent{\bf Acknowledgments.} The author would like to
thank Professor Victor J. W. Guo for helpful comments on this paper.

\end{document}